\documentclass{ifacconf}
\usepackage{graphicx}      
\usepackage{natbib}        
\usepackage{amsmath}
\usepackage{amsfonts}
\usepackage{float}
\usepackage{listings}
\def\postbreak{\raisebox{0ex}[0ex][0ex]{\ensuremath{\hookrightarrow\space}}}
\begin{document}
\begin{frontmatter}
\title{Nonlinear analysis of PLL by the harmonic balance method.}

\author[spb]{Kudryashova E.~V.}
\author[spb,fin]{Kuznetsov N.~V.}
\author[spb,ipmash]{Leonov G.~A.}
\author[spb]{Yuldashev M.~V.}
\author[spb]{Yuldashev R.~V.}

\address[spb]{Faculty of Mathematics and Mechanics,\
Saint-Petersburg State University, Russia}
\address[fin]{Dept. of Mathematical Information Technology,\
  University of Jyv\"{a}skyl\"{a}, Finland} 
\address[ipmash]{Institute of Problems of Mechanical Engineering RAS, Russia}

\begin{abstract}                
  In this paper we discuss the application of
the harmonic balance method
for the global analysis of the classical phase-locked loop (PLL) circuit.
The harmonic balance is non rigorous method,
which is widely used 
for the computation of periodic solutions and the checking of global stability.
The proof of the absence of periodic solutions
is a key step to establish
the global stability of PLL and estimate the pull-in range
(which is an interval of the frequency deviations
such that any solution tends to one of the equilibria).
The advantages and limitations
of the study of the classical PLL with lead-lag filter
using the harmonic balance method is discussed.
\end{abstract}

\end{frontmatter}
\section{Introduction}
Phase-locked loop (PLL) is a nonlinear control system,
which various modifications are widely used in telecommunication and computer architecture
for the master-slave synchronization of oscillators and data demodulations.
Rigorous analysis of the mathematical models of PLLs is
a challenging task and, thus, the simulation and non rigorous methods are often
used in engineering literature for their analysis.

In this paper we discuss the application of
the harmonic balance (HB) method
for the global analysis of the classical PLL.
The harmonic balance is a non-rigorous analytical method,
which allows to study periodic solutions in control systems.
It is widely applied for the study of PLL
(see, e.g. \cite{Margaris-2004,Suarez-2012,Razavi-hb-2016}).
The proof of the absence of periodic solutions
is a key step to establish the global stability
of the PLL model and estimate the pull-in range
(which is an interval of the frequency deviations
such that any solution tends to one of the equilibria).
It is known that the harmonic balance method may lead to wrong conclusion on the global stability,
e.g. it states that well-known Aizerman's and Kalman's conjectures
on the global stability of nonlinear control systems
are valid, while there are known counterexamples
with \emph{hidden oscillations}
(see, e.g. \cite{Pliss-1958,Fitts-1966,Barabanov-1988,BernatL-1996,LeonovBK-2010-DAN,BraginVKL-2011,LeonovK-2011-DAN,LeonovK-2013-IJBC};
the corresponding discrete examples are considered in \cite{Alli-Oke-2012-cu,HeathCS-2015}).
Below we consider advantages and limitations
of the study of classical PLL with lead-lag filter
using the harmonic balance method.
Section~2 introduces the mathematical model of PLL
in a signal's phase space \citep{LeonovKYY-2012-TCASII,LeonovKYY-2015-TCAS}.
In Section~3 the harmonic balance equations are derived,
in Section~4 the harmonic balance equations
are solved numerically
and the obtained results are compared with
the result of the direct simulation of the model.

\section{Classical nonlinear mathematical models
of PLL-based circuits in a signal's phase space}

\begin{figure}[H]
 \centering
 \includegraphics[width=0.73\linewidth]{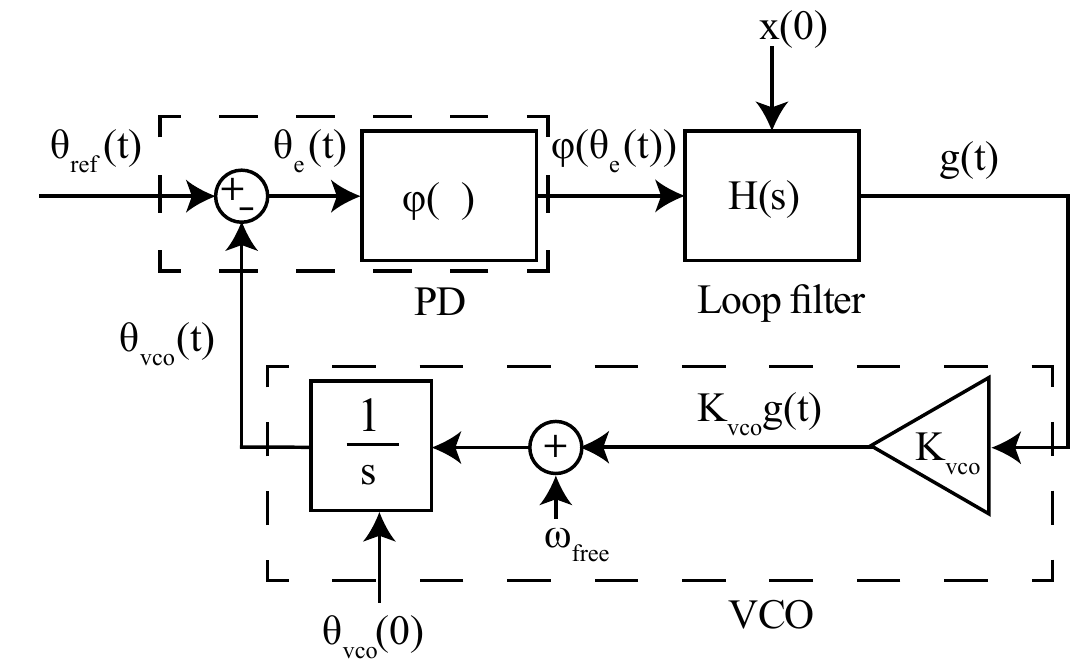}
 \caption{
   PLL-based circuit in a signal's phase space.
 }
 \label{costas-pll-sim-model}
\end{figure}

In classical engineering works (see \citep{Viterbi-1966,Gardner-1966,Best-2007})
various analog PLL-based circuits
are represented in a \emph{signal's phase space} \citep{LeonovKYY-2015-TCAS}
(also named \emph{frequency-domain} \cite[p.338]{Davis-2011})
by a block diagram shown in Fig.~\ref{costas-pll-sim-model}.
Here the Phase Detector (PD) is a nonlinear block;
the phases $\theta_{\rm ref,vco}(t)$ of
the input (reference) and voltage controlled oscillator (VCO) signals
are the PD block inputs,
and the output is the function
$\varphi(\theta_e(t)) = \varphi(\theta_{\rm ref}(t)-\theta_{\rm vco}(t))$
called a \emph{phase detector characteristic},
where
\begin{equation}
  \label{theta_delta_def}
  \begin{aligned}
    & \theta_e(t) = \theta_{\rm ref}(t) - \theta_{\rm vco}(t),
  \end{aligned}
\end{equation}
is called the phase error.
For the classical PLL-based circuits with sinusoidal signal's
waveforms the phase-detector characteristics is sinusoidal:
\begin{equation}\label{phisin}
\begin{aligned}
   & \varphi(\theta_e) = \frac{1}{2}\sin(\theta_e).
\end{aligned}
\end{equation}
The relationship between the input $\varphi(\theta_e(t))$
and the output $g(t)$ of the linear filter (Loop filter) is as follows:
\begin{equation}\label{loop-filter}
 \begin{aligned}
 & \dot x = A x + b \varphi(\theta_e(t)),
 \ g(t) = c^*x + h\varphi(\theta_e(t)),
 \end{aligned}
\end{equation}
where $A$ is a constant $n$-by-$n$ matrix,
$x(t) \in \mathbb{R}^n$ is the filter state,
$x(0)$ is the initial state of filter,
$b$ and $c$ are constant vectors, and $h$ is a number.
The filter transfer function is:\footnote{
  In the control theory the transfer function is often defined
  with the opposite sign
  (see, e.g. \citep{LeonovKYY-2015-TCAS}): $H(s) = c^*(A-sI)^{-1}b-h.$
}
\begin{equation}
\label{transfer-function}
  H(s) = -c^*(A-sI)^{-1}b+h.
\end{equation}
A lead-lag filter \citep{Best-2007} (with $H(0) = 1$),
or a PI filter ($H(0)$ is infinite)
is usually used as the loop filter.
The control signal $g(t)$ adjusts the VCO frequency to
the frequency of the input signal:
\begin{equation} \label{vco first}
   \dot\theta_{\rm vco}(t) = \omega_{\rm vco}(t) = \omega_{\rm vco}^{\text{free}}
   + K_{\rm vco}g(t),
\end{equation}
where $\omega_{\rm vco}^{\text{free}}$ is the VCO free-running frequency
(i.e. for $g(t)\equiv 0$) and $K_{\rm vco}$ is the VCO gain.
Nonlinear VCO models can be considered similarly, see,
e.g. \citep{Margaris-2004,BianchiKLYY-2016}.
The frequency of the input signal (reference frequency) is usually assumed
to be constant:
\begin{equation}\label{omega1-const}
  \dot\theta_{\rm ref}(t) = \omega_{\rm ref}(t) \equiv \omega_{\rm ref}.
\end{equation}
The difference between the reference frequency and the VCO free-running
frequency is denoted as $\omega_e^{\text{free}}$:
\begin{equation}
  \label{omega_delta_def}
  \begin{aligned}
    & \omega_e^{\text{free}} \equiv \omega_{\rm ref} - \omega_{\rm vco}^{\text{free}}.
  \end{aligned}
\end{equation}
Combining equations \eqref{theta_delta_def}, \eqref{loop-filter},
and \eqref{vco first}--\eqref{omega_delta_def}, we get
\begin{equation}\label{system1}
 \begin{aligned}
   & \dot\theta_e = \omega_e^{\text{free}}
   - K_{\rm vco} g(t).
 \end{aligned}
\end{equation}

By \eqref{loop-filter} and \eqref{system1} we obtain
a \emph{nonlinear mathematical model in a signal's phase space}
(i.e. in the state space: the filter's state $x$
and the difference between the signal's phases $\theta_e$):
\begin{equation}\label{final_system}
 \begin{aligned}
 & \dot x = A x + b \varphi(\theta_e(t)),
 \\
   & \dot\theta_e = \omega_e^{\text{free}}
   - K_{\rm vco} \big(c^*x + h\varphi(\theta_e(t))\big).
 \end{aligned}
\end{equation}

In the case of PD characteristic \eqref{phisin},
system  \eqref{final_system} is not changed under the transformation
\begin{equation}\label{odd-change}
  \big(\omega_e^{\text{free}},x(t),\theta_e(t)) \rightarrow
  \big(-\omega_e^{\text{free}},-x(t),-\theta_e(t))
\end{equation}
and, thus, we can analyze system \eqref{final_system}
only with $\omega_e^{\text{free}}>0$
and introduce the concept of \emph{frequency deviation}
\begin{center}{$|\omega_e^{\text{free}}| =
  |\omega_{\rm ref} - \omega_{\rm vco}^{\text{free}}|$.}
\end{center}

The pull-in range is a widely used engineering concept
(see, e.g. \citep[p.40]{Gardner-1966}, \citep[p.61]{Best-2007}).
The following rigorous definition is suggested \citep{KuznetsovLYY-2015-IFAC-Ranges,LeonovKYY-2015-TCAS,BestKLYY-2016}.
The largest interval of frequency deviations
$0 \leq |\omega_{e}^{\rm{free}}|<\omega_{\rm{pull-in}}$
such that the nonlinear mathematical model of PLL in the signal's phase space
acquires lock for arbitrary initial phase difference and filter state
(i.e. any trajectory tends to an equilibrium point)
is called a \emph{pull-in range}, $\omega_{\rm{pull-in}}$ is called
a \emph{pull-in frequency}.

This definition implies that for any frequency deviation
from pull-in range
the mathematical model of PLL does not contain periodic solutions.
This property can be used to obtain necessary conditions of pull-in range
(see, e.g. \citep{Razavi-hb-2016,BianchiKLYY-2015,BianchiKLYY-2016}).
In the next section the application of harmonic balance
method to the PLL with lead-lag filter is discussed.

\section{Harmonic Balance method}
\begin{figure*}
     \includegraphics[scale=0.3]{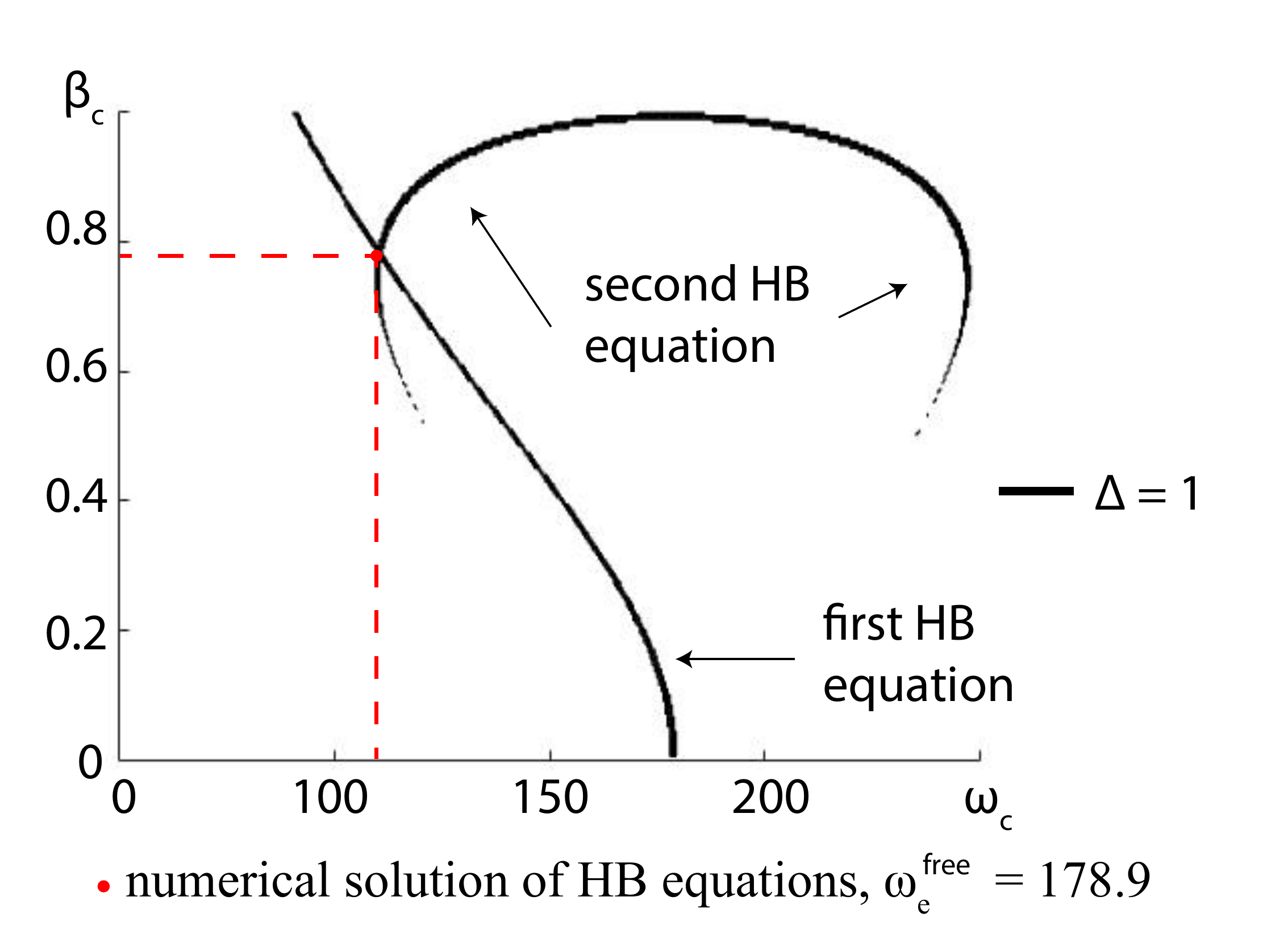}
     \includegraphics[scale=0.3]{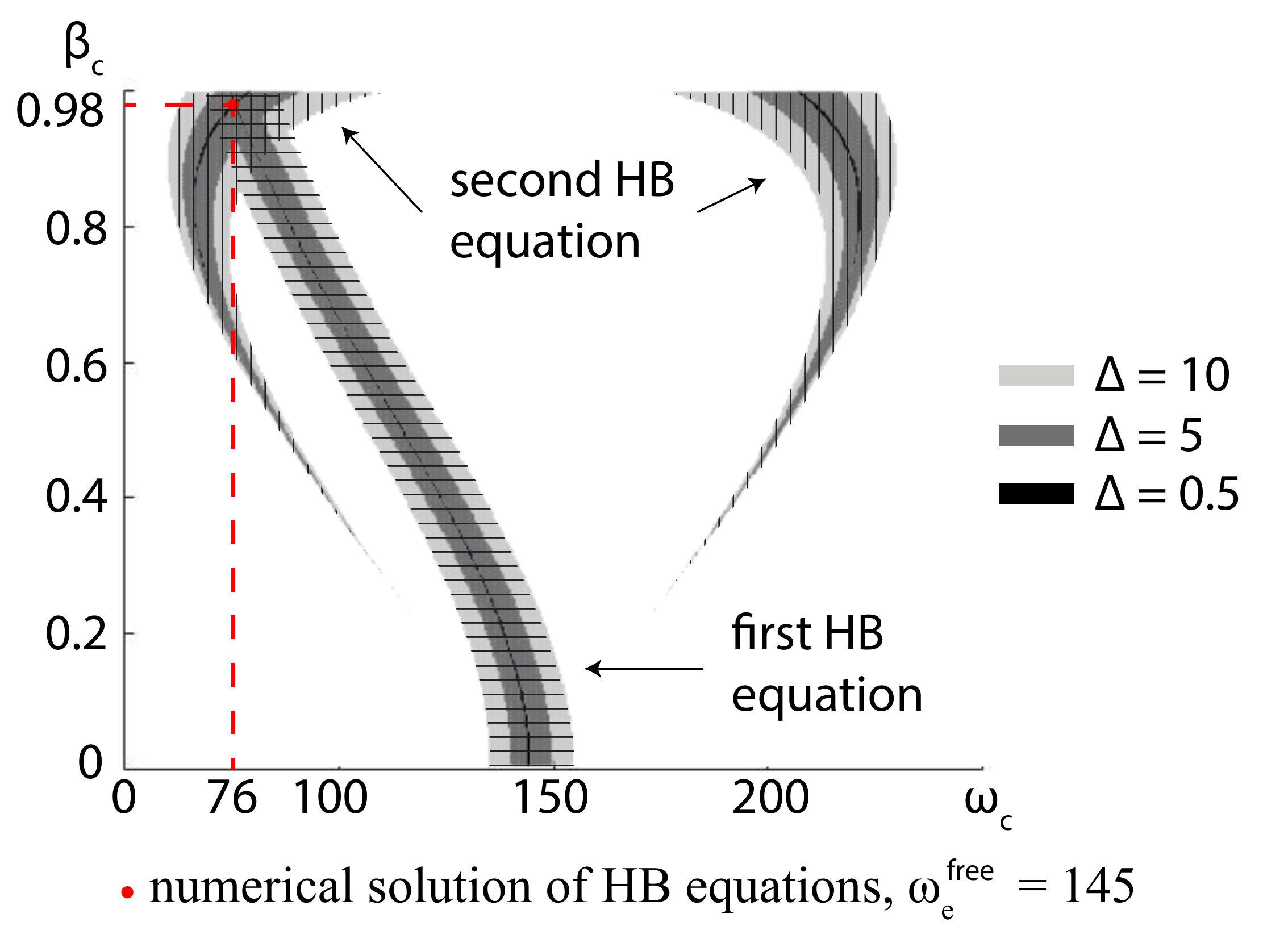}
     \centering
     \caption{Region of $\omega_c$ and $\beta_c$,
     where HB equations \eqref{HB-equations-simpl} are satisfied with tolerance $\Delta$.
     $\Delta$ is a maximum absolute difference between the right-hand
    side and left-hand-side of corresponding equations.
  The solution of HB equations for $\omega_e^{\rm free} = 178.9$ is as follows
                $\omega_c \approx 110$,
                $\beta_c \approx 0.785$,
                $\theta_c \approx 0.7088$.
  The solution of HB equations for $\omega_e^{\rm free} = 145$ has the form
                $\omega_c \approx 76$,
                $\beta_c \approx 0.98$,
                $\theta_c \approx 0.88$.
                Parameters:
                        $K_{\rm vco} = 250$,
                        $\tau_1 = 0.0448$,
                        $\tau_2 = 0.0185$.}
     \label{omega 178}
 \end{figure*}

Following \citep{Shahgildyan-1972},
let us look for a solution in the following form
\begin{equation} \label{solution}
\begin{aligned}
  & \theta_e(t) = \omega_{c}t + \theta_c + \frac{\pi}{2}  - \beta_c \sin(\omega_{c}t).
\end{aligned}
\end{equation}
Here $\omega_{c}$, $\theta_c$, and $\beta_c$ are unknown parameters of the solution.
The output of sinusoidal phase detector has the form
\begin{equation}
\begin{aligned}
& \varphi(\theta_e(t)) = \frac{1}{2}\sin(\omega_{c}t
+ \theta_c + \frac{\pi}{2} - \beta_c \sin(\omega_{c}t))
\\
& = \frac{1}{2}\cos(\omega_{c}t + \theta_c - \beta_c \sin(\omega_{c}t)).
\end{aligned}
\end{equation}
By using the first two elements of the following equations from Bessel functions
theory  (Jacobi-Anger expansion \citep{abramowitz1964handbook}):
\begin{equation}
\begin{aligned}
& \cos(\beta_c \sin(\omega_{c}t)) =
    J_0(\beta_c)
    + 2\sum\limits_{n=1}^{\infty}J_{2n}(\beta_c)\cos(2n\omega_{c}t),
    \\
& \sin(\beta_c \sin(\omega_{c}t)) =
    2\sum\limits_{n=0}^{\infty}J_{2n+1}(\beta_c)\sin((2n+1)\omega_{c}t),
    \\
\end{aligned}
\end{equation}
we obtain the following approximation
\begin{equation}
\begin{aligned}
& \cos(\omega_{c}t + \theta_c - \beta_c \sin(\omega_{c}t))
\\
\approx
\\
& \Big(\cos(\omega_{c}t)\cos(\theta_c) - \sin(\omega_{c}t)\sin(\theta_c)\Big)
\\
& \Big(J_0(\beta_c)
    + 2J_{2}(\beta_c)\cos(2\omega_{c}t)
    \Big) +
\\
&  + \Big(\sin(\omega_{c}t)\cos(\theta_c) + \cos(\omega_{c}t)\sin(\theta_c)\Big)
\\
& 2J_{1}(\beta_c)\sin(\omega_{c}t)
\\
=
\\
& J_0(\beta_c)\cos(\omega_{c}t)\cos(\theta_c)
\\
& - J_0(\beta_c)\sin(\omega_{c}t)\sin(\theta_c)
\\
&    + J_{2}(\beta_c)\Big(
            \cos(\omega_{c}t)
            +
            \cos(3\omega_{c}t)
        \Big)\cos(\theta_c)
\\
&
     - J_{2}(\beta_c)\Big(
            \sin(-\omega_{c}t)
            +
            \sin(3\omega_{c}t)
        \Big)
     \sin(\theta_c)
\\
& + J_{1}(\beta_c)\Big(1-\cos(2\omega_{c}t)\Big)\cos(\theta_c)
\\
&  + J_{1}(\beta_c)\sin(2\omega_{c}t)\sin(\theta_c).
\end{aligned}
\end{equation}
Excluding higher harmonics ($\cos(k\omega_{c}t)$
and $\sin(k\omega_{c}t)$ for $k \geq 2$)
we get:
\begin{equation}
\begin{aligned}
& \cos(\omega_{c}t + \theta_c - \beta_c \sin(\omega_{c}t))
\\
\approx &
\\
& J_0(\beta_c)\cos(\omega_{c}t)\cos(\theta_c)
\\
& - J_0(\beta_c)\sin(\omega_{c}t)\sin(\theta_c)
\\
&
     +  J_{2}(\beta_c)
            \cos(\omega_{c}t)
     \cos(\theta_c)
\\
&    + J_{2}(\beta_c)
            \sin(\omega_{c}t)
        \sin(\theta_c)
\\
&  + J_{1}(\beta_c)\cos(\theta_c)
\\
=
\\
& \Big(J_0(\beta_c) + J_{2}(\beta_c)\Big)\cos(\omega_{c}t)\cos(\theta_c)
\\
& - \Big(J_{0}(\beta_c) - J_2(\beta_c)\Big)\sin(\omega_{c}t)\sin(\theta_c)
\\
&  + J_{1}(\beta_c)\cos(\theta_c).
\end{aligned}
\end{equation}
Then the output of the linear loop filter can be approximated
 as\footnote{Here filter \eqref{transfer-function}
is considered as a linear time-invariant (LTI) system}
\begin{equation}
\label{filter approx}
\begin{aligned}
& g(t) \approx
\frac{1}{2}J_1(\beta_c)\cos(\theta_c)
\\
&
- \frac{1}{2}|H(i\omega_c)|
    [J_0(\beta_c) - J_2(\beta_c)]
    \sin(\theta_c)\sin(\omega_{c}t - \psi_{\omega})
\\
&
+ \frac{1}{2}|H(i\omega_c)|
    [J_0(\beta_c) + J_2(\beta_c)]
    \cos(\theta_c)\cos(\omega_{c}t - \psi_{\omega}),
\end{aligned}
\end{equation}
where $\psi_{\omega}$ and $|H(i\omega_c)|$
are the filter phase shift and gain for the frequency $\omega_c$.
Taking derivative of \eqref{solution}, we get
\begin{equation}
\label{solution derivative}
\begin{aligned}
& \dot\theta_e(t) = \omega_{c} - \beta_c \omega_c \cos(\omega_{c}t).
\end{aligned}
\end{equation}
Substituting \eqref{filter approx}
and \eqref{solution derivative} in PLL equation \eqref{system1}, we have
\begin{equation}
\label{almost hb}
\begin{aligned}
& \omega_{c}- \omega_{c}\beta_c \Big(
        \cos(\omega_{c}t - \psi_{\omega})\cos(\psi_{\omega})
        -
        \sin(\omega_{c}t - \psi_{\omega})\sin(\psi_{\omega})
        \Big)
\\
& =
 \omega_e^{\rm free} - \frac{K_{\rm vco}}{2}
 \Big(
J_1(\beta_c)\cos(\theta_c)
\\
&
- |H(i\omega_c)|
    [J_0(\beta_c) - J_2(\beta_c)]
    \sin(\theta_c)\sin(\omega_{c}t - \psi_{\omega})
\\
&
+ |H(i\omega_c)|
    [J_0(\beta_c) + J_2(\beta_c)]
    \cos(\theta_c)\cos(\omega_{c}t - \psi_{\omega})
 \Big).
\end{aligned}
\end{equation}
By equations \eqref{almost hb} we get the following harmonic balance equations
\begin{equation}
\label{hb eq}
\begin{aligned}
& \omega_{c}+ \frac{K_{\rm vco}}{2}J_1(\beta_c)\cos(\theta_c) = \omega_e^{\rm free},
\\
& \omega_{c}\beta_c \cos(\psi_{\omega}) =
    \frac{K_{\rm vco}}{2}|H(i\omega_c)|[J_0(\beta_c) + J_2(\beta_c)]
    \cos(\theta_c),
\\
& \omega_{c}\beta_c \sin(\psi_{\omega}) =
    \frac{K_{\rm vco}}{2}|H(i\omega_c)|[J_0(\beta_c) - J_2(\beta_c)]
    \sin(\theta_c).
\end{aligned}
\end{equation}
Using the property of Bessel functions:
\begin{equation}
\begin{aligned}
& J_0(\beta_c) + J_2(\beta_c) = \frac{2J_1(\beta_c)}{\beta_c},
\end{aligned}
\end{equation}
we have
\begin{equation}
\label{HB equations}
\begin{aligned}
& \cos(\theta_c) =
\frac{\omega_e^{\rm free} - \omega_{c}}{\frac{K_{\rm vco}}{2}J_1(\beta_c)},
\\
& \omega_{c}\beta_c \cos(\psi_{\omega}) =
    \frac{K_{\rm vco}}{2}|H(i\omega_c)|\frac{2J_1(\beta_c)}{\beta_c}
    \frac{\omega_e^{\rm free} - \omega_{c}}{\frac{K_{\rm vco}}{2}J_1(\beta_c)},
\\
& \omega_{c}\beta_c \sin(\psi_{\omega}) =
    \frac{K_{\rm vco}}{2}|H(i\omega_c)|[J_0(\beta_c) - J_2(\beta_c)]
    \cos(\theta_c),
\end{aligned}
\end{equation}
which is equal to the following
\begin{equation}
\label{HB equations}
\begin{aligned}
& \sin(\theta_c) = \frac{\omega_e^{\rm free}
- \omega_{c}}{\frac{K_{\rm vco}}{2}J_1(\beta_c)},
\\
& \omega_{c} =
    \frac{2|H(i\omega_c)|}{\beta_c^2\cos(\psi_{\omega})}
    (\omega_e^{\rm free} - \omega_{c}),
\\
& \omega_{c}\beta_c \sin(\psi_{\omega}) =
    \frac{K_{\rm vco}}{2}|H(i\omega_c)|[J_0(\beta_c) - J_2(\beta_c)]
    \sin(\theta_c).
\end{aligned}
\end{equation}
Finally,
\begin{equation}
\label{HB equations}
\begin{aligned}
& \sin(\theta_c) = \frac{\omega_e^{\rm free}
- \omega_{c}}{\frac{K_{\rm vco}}{2}J_1(\beta_c)},
\\
& \omega_{c} =
    \frac{2|H(i\omega_c)|}{2|H(i\omega_c)| + \beta_c^2\cos(\psi_{\omega})}
    \omega_e^{\rm free},
\\
& \omega_{c}\beta_c \sin(\psi_{\omega}) =
    \frac{K_{\rm vco}}{2}|H(i\omega_c)|[J_0(\beta_c) - J_2(\beta_c)]
    \sin(\theta_c).
\end{aligned}
\end{equation}
Here $J_{0,1,2}$ are Bessel functions;
$\omega_c$, $\beta_c$, and $\theta_c$ are unknown parameters of the solution.

In the next section we consider numerical solution of \eqref{HB equations}.

\section{Numerical solutions of Harmonic-Balance equations
for lead-lag filter} 
\label{sec:numerical_solution_of_hb}
For lead-lag filter we have
\begin{equation}
\begin{aligned}
& |H(i \omega_c)| = \Big|\frac{1+i\tau_2 \omega_c}{1+i\tau_1 \omega_c}\Big|,
\\
& \cos(\psi_{\omega}) = \cos\left(arg\left(\frac{1+i\tau_2 \omega_c}{1+i\tau_1 \omega_c}\right)\right),
\\
& \sin(\psi_{\omega}) = \sin\left(arg\left(\frac{1+i\tau_2 \omega_c}{1+i\tau_1 \omega_c}\right)\right).
\end{aligned}
\end{equation}

 Let us find numerically the solutions of \eqref{HB equations}.
 To solve nonlinear equations \eqref{HB equations},
 it is possible to apply MATLAB function ``vpasolve'',
 but the result depends on the \emph{initial guess}
 that is not convenient for the checking of
 the absence of solutions of \eqref{HB equations}.
 Thus, we consider the difference between
 the right-hand side and left-hand side of equations \eqref{HB equations}.
 Since we cannot find exact solution,
 we plot the points on $(\omega_c, \beta_c)$-plane for which
 the absolute value of the differences between the right-hand side and left-hand-side of
\eqref{HB equations} is less than $\Delta$, i.e.
 \begin{equation}
\label{HB-equations-simpl}
\begin{aligned}
& \Big|
    \omega_{c} - \omega_e^{\rm free}
    \frac{2|H(i\omega_c)|}{\beta_c^2 \cos(\psi_{\omega}) + 2|H(i\omega_c)|}
   \Big|
   < \Delta,
\\
& \Big|
    \omega_{c}\beta_c \sin(\psi_{\omega}) -
    K_{\rm vco}|H(i\omega_c)|\big(J_0(\beta_c) - J_2(\beta_c)\big)
    \\
    & \quad\quad\quad\quad\quad\quad\quad\quad\quad\quad
    \sqrt{1 - \Big(\frac{\omega_e^{\rm free} - \omega_{c}}{K_{\rm vco}J_1(\beta_c)}\Big)^2}
    \Big| < \Delta.
\\
\end{aligned}
\end{equation}
The values $(\omega_c, \beta_c)$ satisfying
conditions \eqref{HB-equations-simpl}
with $\Delta = 1$ are shown in Fig.~\ref{omega 178}.
In the left subfigure there are two areas,
which corresponds to the first and the second equations of \eqref{HB-equations-simpl}.
The intersection of this areas gives an approximation of solution
of harmonic-balance equations, e.g. the intersection contains
the following point $\omega_c \approx 76$, $\beta_c \approx 0.98$,
$\theta_c \approx 0.88$.
For the parameters $\omega_c \approx 76$, $\beta_c \approx 0.98$,
$\theta_c \approx 0.88$ we plot \eqref{solution}
that is the results of numerical simulation of system \eqref{final_system}
with zero initial conditions.
  \begin{figure}[h]
     \includegraphics[width=0.4\textwidth]{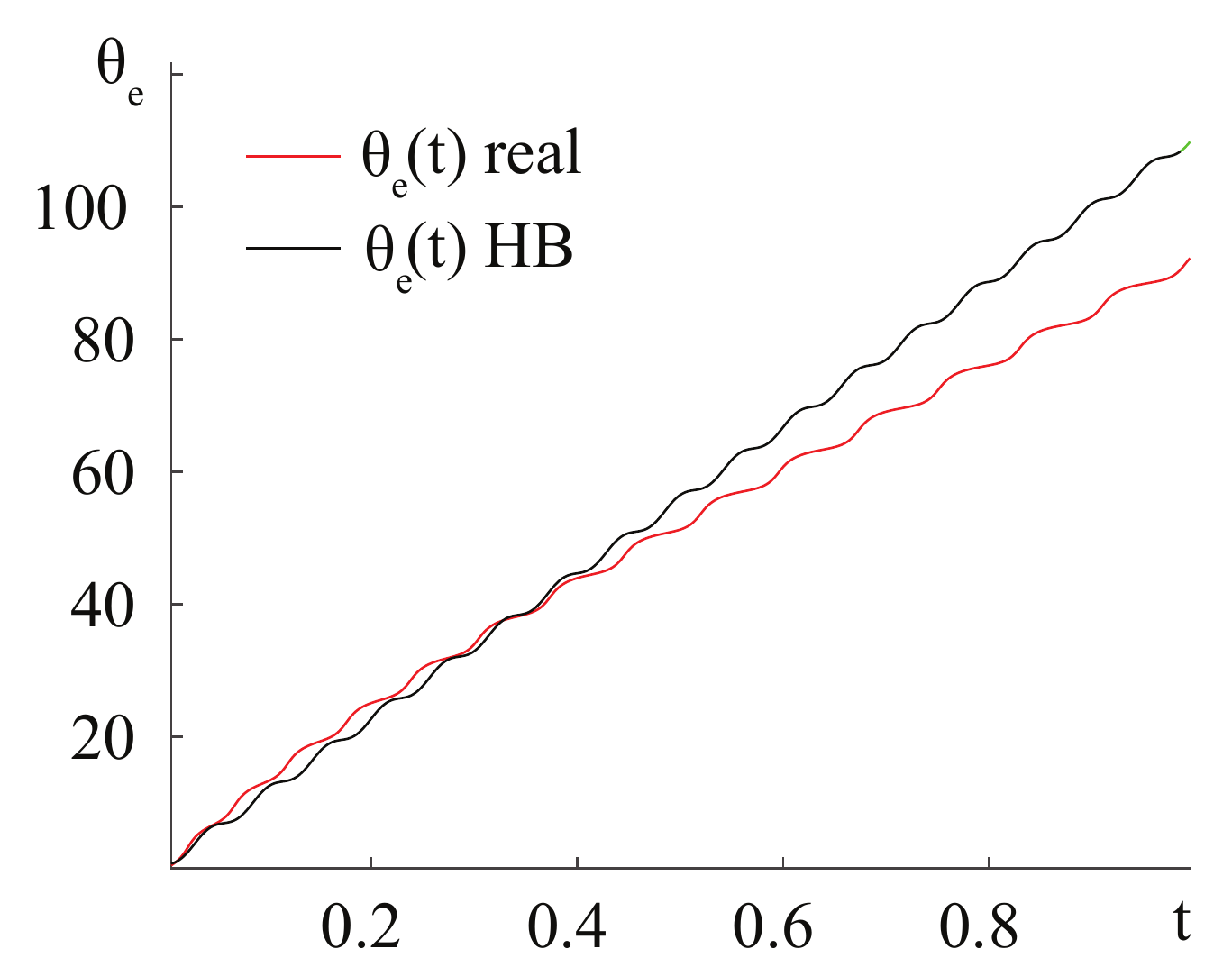}
     \centering
     \caption{Solution obtained by HB vs real solution.
     Solution of HB equations:
                $\omega_c \approx 110$,
                $\beta_c \approx 0.785$,
                $\theta_c \approx 0.7088$.
                Parameters:
                        $\omega_e^{\rm free} = 178.9$,
                        $K_{\rm vco} = 250$,
                        $\tau_1 = 0.0448$,
                        $\tau_2 = 0.0185$.}
     \label{omega 178 theta}
 \end{figure}
As shown in Fig.~\ref{omega 178 theta} the solution $\theta_e(t)$ tends
to infinity and the approximation given by the harmonic balance method
is correct and contains periodic part (cycle).

If we consider smaller values of $\omega_e^{\rm free}$
(up to $\omega_e^{\rm free} = 145$), then equations \eqref{HB-equations-simpl}
still have a solution (see right subfigure in Fig.~\ref{omega 178}).
However in this case the harmonic balance method leads us to a wrong conclusion
since we can not reveal corresponding cycle in system \eqref{final_system}
by direct simulation
(see the comparison of numerical solutions in Fig.~\ref{omega 145 theta})).
 \begin{figure}[h]
     \includegraphics[width=0.4\textwidth]{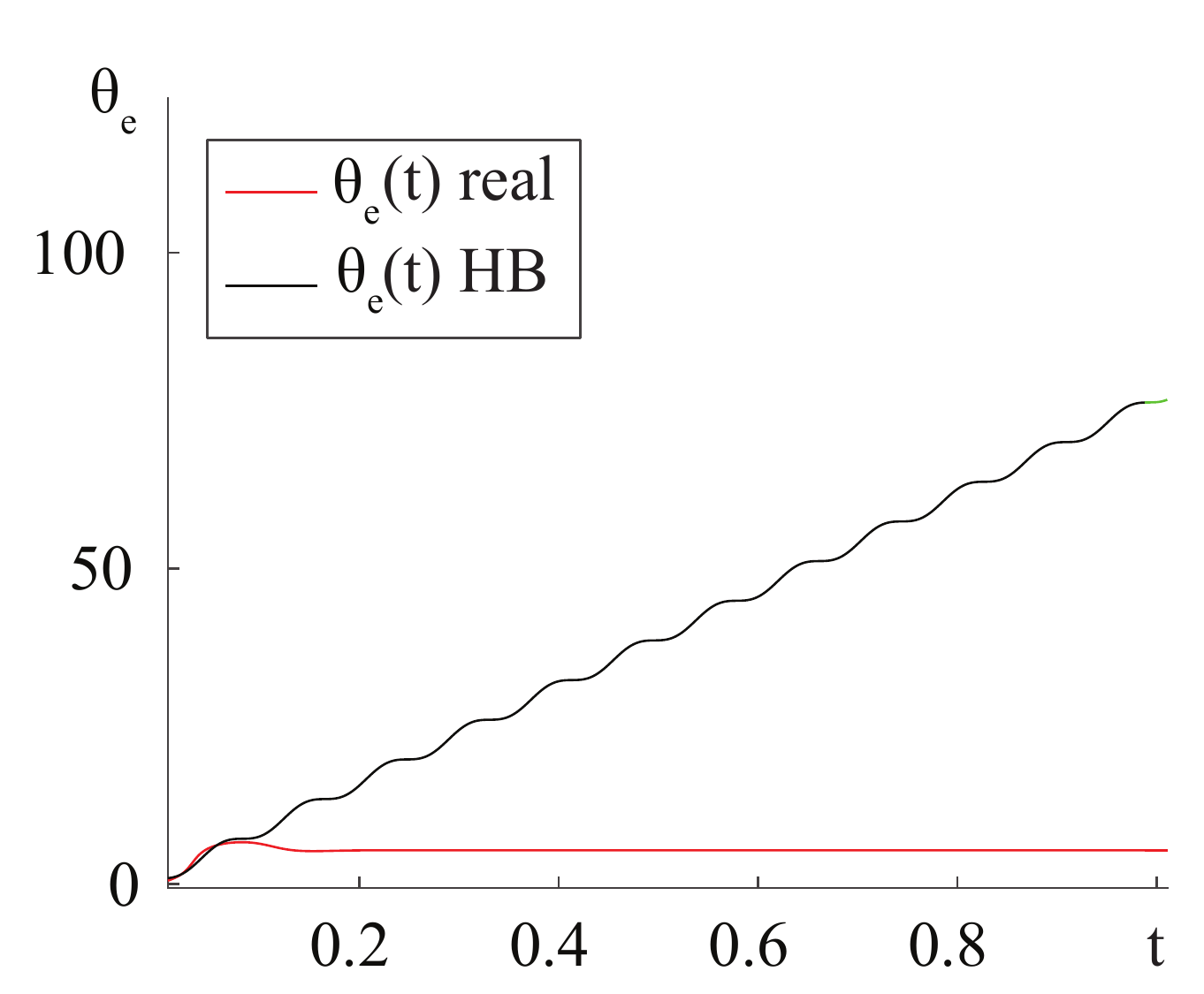}
     \centering
     \caption{Solution obtained by HB vs real solution.
     The solution of HB equations:
                $\omega_c \approx 76$,
                $\beta_c \approx 0.98$,
                $\theta_c \approx 0.88$.
                Parameters:
                        $\omega_e^{\rm free} = 145$,
                        $K_{\rm vco} = 250$,
                        $\tau_1 = 0.0448$,
                        $\tau_2 = 0.0185$.}
      \label{omega 145 theta}
 \end{figure}

Also it is possible to check that harmonic balance equations \eqref{HB equations}
have a solution for any parameters.
But the solutions with $\beta_c > 1$ is usually excluded
(\citep{ShahgildyanL-1966,Shahgildyan-1972})
because the phase is supposed to be nonnegative.
If there exists a solution for $0 < \beta < 1$,
then HB implies the existence of cycle \eqref{solution}.
The frequency of the cycle is limited by a cut-off frequency of
the filter $0 < \omega_c < \omega_{\text{cut-off}}$ and $\theta_c \in [0,2\pi]$.

Remark that simulation itself may not reveal a
complex behavior of PLL: such examples, where the simulation of PLL-based circuits leads to unreliable results, are demonstrated in \citep{BianchiKLYY-2015,BlagovKLYY-2015,KuznetsovLYY-2017-CNSNS}.
Consider $\omega_e^{\rm free} = 178.545$  and the lead-lag filter with $\tau_1 = 0.0448$, $\tau_2 = 0.0185$.
This value is close to bifurcation point,
where a periodic oscillations appears.
Simulation with relatively small precision ('MaxStep', $0.01$, 'RelTol', $2e-3$, 'AbsTol', $2e-3$) shows absence of cycles,
while simulation with precision ('MaxStep', $0.001$, 'RelTol', $2e-6$, 'AbsTol', $2e-6$)
allows to reveal a cycle (see Fig.~\ref{omega 145 theta}).
\begin{figure}[H]
     \includegraphics[width=0.4\textwidth]{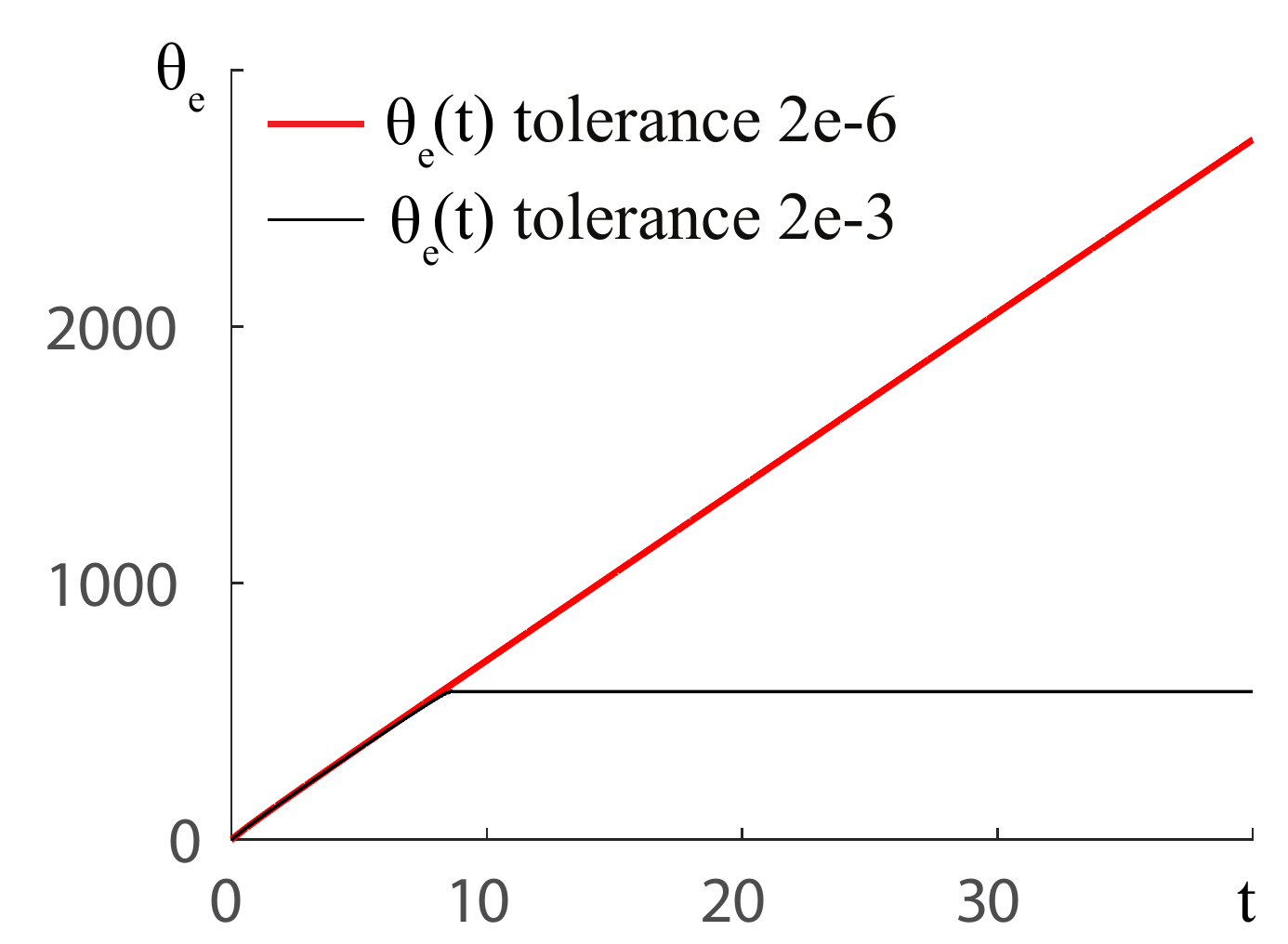}
     \centering
     \caption{
     $\omega_e^{\rm free} = 178.545$,
                        $K_{\rm vco} = 250$,
                        $\tau_1 = 0.0448$,
                        $\tau_2 = 0.0185$.
     MATLAB `odeset' parameters:
     black line --- odeset('MaxStep', $0.01$, 'RelTol', $2e-3$, 'AbsTol', $2e-3$),
     red (grey) line --- odeset('MaxStep', $0.001$, 'RelTol', $2e-6$, 'AbsTol', $2e-6$).}
      \label{omega 145 theta}
 \end{figure}

This example demonstrates the difficulties of
numerical search of so-called hidden oscillations,
whose basin of attraction does not overlap with the neighborhood
of the equilibrium point, and thus it may be difficult to find them numerically
\citep{LeonovK-2013-IJBC,LeonovKM-2015-EPJST,Kuznetsov-2016}.
In this case the observation of one or another
stable solution may depend on the initial data and integration step.

\section{Conclusions}
While harmonic balance method is widely used for the estimation of the pull-in range,
it may lead to wrong results.
Corresponding examples are discussed in the paper.
The pull-in range of PLL-based circuits with first-order filters
can be estimated
using phase plane analysis methods \citep{Tricomi-1933,AndronovVKh-1937,Shakhtarin-1969,Belyustina-1970-eng}.
For rigorous nonlinear analysis of multidimensional PLL models
one may use special modifications of the classical stability criteria developed for the cylindrical phase space in \citep{LeonovKYY-2015-TCAS,LeonovK-2014-book}.

\section*{Acknowledgment}
This work was supported by Russian Science Foundation (project 14-21-00041).

\bibliographystyle{ifacconf}

\begin{thebibliography}{35}
\providecommand{\natexlab}[1]{#1}
\providecommand{\url}[1]{\texttt{#1}}
\providecommand{\urlprefix}{URL }
\expandafter\ifx\csname urlstyle\endcsname\relax
  \providecommand{\doi}[1]{doi:\discretionary{}{}{}#1}\else
  \providecommand{\doi}{doi:\discretionary{}{}{}\begingroup
  \urlstyle{rm}\Url}\fi

\bibitem[{Abramowitz and Stegun(1964)}]{abramowitz1964handbook}
Abramowitz, M. and Stegun, I.A. (1964).
\newblock \emph{Handbook of mathematical functions: with formulas, graphs, and
  mathematical tables}, volume~55.
\newblock Courier Corporation.

\bibitem[{Alli-Oke et~al.(2012)Alli-Oke, Carrasco, Heath, and
  Lanzon}]{Alli-Oke-2012-cu}
Alli-Oke, R., Carrasco, J., Heath, W., and Lanzon, A. (2012).
\newblock A robust {K}alman conjecture for first-order plants.
\newblock \emph{IFAC Proceedings Volumes (IFAC-PapersOnline)}, 7, 27--32.
\newblock \doi{10.3182/20120620-3-DK-2025.00161}.

\bibitem[{Andronov et~al.(1937)Andronov, Vitt, and Khaikin}]{AndronovVKh-1937}
Andronov, A.A., Vitt, E.A., and Khaikin, S.E. (1937).
\newblock \emph{Theory of Oscillators (in Russian)}.
\newblock ONTI NKTP SSSR.
\newblock [English transl.: 1966, Pergamon Press].

\bibitem[{Barabanov(1988)}]{Barabanov-1988}
Barabanov, N.E. (1988).
\newblock On the {K}alman problem.
\newblock \emph{Sib. Math. J.}, 29(3), 333--341.

\bibitem[{Belyustina et~al.(1970)Belyustina, Brykov, Kiveleva, and
  Shalfeev}]{Belyustina-1970-eng}
Belyustina, L., Brykov, V., Kiveleva, K., and Shalfeev, V. (1970).
\newblock On the magnitude of the locking band of a phase-shift automatic
  frequency control system with a proportionally integrating filter.
\newblock \emph{Radiophysics and Quantum Electronics}, 13(4), 437--440.

\bibitem[{Bernat and Llibre(1996)}]{BernatL-1996}
Bernat, J. and Llibre, J. (1996).
\newblock Counterexample to {K}alman and {M}arkus-{Y}amabe conjectures in
  dimension larger than 3.
\newblock \emph{Dynamics of Continuous, Discrete and Impulsive Systems}, 2(3),
  337--379.

\bibitem[{Best(2007)}]{Best-2007}
Best, R. (2007).
\newblock \emph{Phase-Locked Loops: Design, Simulation and Application}.
\newblock McGraw-Hill, 6th edition.

\bibitem[{Best et~al.(2016)Best, Kuznetsov, Leonov, Yuldashev, and
  Yuldashev}]{BestKLYY-2016}
Best, R., Kuznetsov, N., Leonov, G., Yuldashev, M., and Yuldashev, R. (2016).
\newblock Tutorial on dynamic analysis of the {C}ostas loop.
\newblock \emph{Annual Reviews in Control}, 42, 27--49.
\newblock \doi{10.1016/j.arcontrol.2016.08.003}.

\bibitem[{Bianchi et~al.(2016{\natexlab{a}})Bianchi, Kuznetsov, Leonov,
  Seledzhi, Yuldashev, and Yuldashev}]{BianchiKLYY-2016}
Bianchi, G., Kuznetsov, N., Leonov, G., Seledzhi, S., Yuldashev, M., and
  Yuldashev, R. (2016{\natexlab{a}}).
\newblock Hidden oscillations in {SPICE} simulation of two-phase {C}ostas loop
  with non-linear {VCO}.
\newblock \emph{IFAC-PapersOnLine}, 49(14), 45--50.
\newblock \doi{10.1016/j.ifacol.2016.07.973}.

\bibitem[{Bianchi et~al.(2016{\natexlab{b}})Bianchi, Kuznetsov, Leonov,
  Yuldashev, and Yuldashev}]{BianchiKLYY-2015}
Bianchi, G., Kuznetsov, N., Leonov, G., Yuldashev, M., and Yuldashev, R.
  (2016{\natexlab{b}}).
\newblock Limitations of {PLL} simulation: hidden oscillations in {MATLAB} and
  {SPICE}.
\newblock \emph{International Congress on Ultra Modern Telecommunications and
  Control Systems and Workshops (ICUMT 2015)}, 2016-January, 79--84.
\newblock \doi{10.1109/ICUMT.2015.7382409}.

\bibitem[{Blagov et~al.(2016)Blagov, Kuznetsov, Leonov, Yuldashev, and
  Yuldashev}]{BlagovKLYY-2015}
Blagov, M., Kuznetsov, N., Leonov, G., Yuldashev, M., and Yuldashev, R. (2016).
\newblock Simulation of {PLL} with impulse signals in {MATLAB}: Limitations,
  hidden oscillations, and pull-in range.
\newblock \emph{International Congress on Ultra Modern Telecommunications and
  Control Systems and Workshops (ICUMT 2015)}, 2016-January, 85--90.
\newblock \doi{10.1109/ICUMT.2015.7382410}.

\bibitem[{Bragin et~al.(2011)Bragin, Vagaitsev, Kuznetsov, and
  Leonov}]{BraginVKL-2011}
Bragin, V., Vagaitsev, V., Kuznetsov, N., and Leonov, G. (2011).
\newblock Algorithms for finding hidden oscillations in nonlinear systems.
  {T}he {A}izerman and {K}alman conjectures and {C}hua's circuits.
\newblock \emph{Journal of Computer and Systems Sciences International}, 50(4),
  511--543.
\newblock \doi{10.1134/S106423071104006X}.

\bibitem[{Davis(2011)}]{Davis-2011}
Davis, W. (2011).
\newblock \emph{Radio Frequency Circuit Design}.
\newblock Wiley Series in Microwave and Optical Engineering. Wiley, IEEE Press.

\bibitem[{Fitts(1966)}]{Fitts-1966}
Fitts, R.E. (1966).
\newblock Two counterexamples to {A}izerman's conjecture.
\newblock \emph{Trans. IEEE}, AC-11(3), 553--556.

\bibitem[{Gardner(1966)}]{Gardner-1966}
Gardner, F. (1966).
\newblock \emph{Phaselock techniques}.
\newblock John Wiley \& Sons, New York.

\bibitem[{Heath et~al.(2015)Heath, Carrasco, and de~la Sen}]{HeathCS-2015}
Heath, W.P., Carrasco, J., and de~la Sen, M. (2015).
\newblock Second-order counterexamples to the discrete-time {K}alman
  conjecture.
\newblock \emph{Automatica}, 60, 140 -- 144.

\bibitem[{Homayoun and Razavi(2016)}]{Razavi-hb-2016}
Homayoun, A. and Razavi, B. (2016).
\newblock On the stability of charge-pump phase-locked loops.
\newblock \emph{IEEE Transactions on Circuits and Systems I: Regular Papers},
  63(6), 741--750.
\newblock \doi{10.1109/TCSI.2016.2537823}.

\bibitem[{Kuznetsov(2016)}]{Kuznetsov-2016}
Kuznetsov, N. (2016).
\newblock Hidden attractors in fundamental problems and engineering models. {A}
  short survey.
\newblock \emph{Lecture Notes in Electrical Engineering}, 371, 13--25.
\newblock \doi{10.1007/978-3-319-27247-4\_2}.
\newblock (Plenary lecture at International Conference on Advanced Engineering
  Theory and Applications 2015).

\bibitem[{Kuznetsov et~al.(2015)Kuznetsov, Leonov, Yuldashev, and
  Yuldashev}]{KuznetsovLYY-2015-IFAC-Ranges}
Kuznetsov, N., Leonov, G., Yuldashev, M., and Yuldashev, R. (2015).
\newblock Rigorous mathematical definitions of the hold-in and pull-in ranges
  for phase-locked loops.
\newblock \emph{IFAC-PapersOnLine}, 48(11), 710--713.
\newblock \doi{10.1016/j.ifacol.2015.09.272}.

\bibitem[{Kuznetsov et~al.(2017)Kuznetsov, Leonov, Yuldashev, and
  Yuldashev}]{KuznetsovLYY-2017-CNSNS}
Kuznetsov, N., Leonov, G., Yuldashev, M., and Yuldashev, R. (2017).
\newblock Hidden attractors in dynamical models of phase-locked loop circuits:
  limitations of simulation in {MATLAB} and {SPICE}.
\newblock \emph{Commun Nonlinear Sci Numer Simulat}, 51, 39--49.
\newblock \doi{10.1016/j.cnsns.2017.03.010}.

\bibitem[{Leonov et~al.(2010)Leonov, Bragin, and Kuznetsov}]{LeonovBK-2010-DAN}
Leonov, G., Bragin, V., and Kuznetsov, N. (2010).
\newblock Algorithm for constructing counterexamples to the {K}alman problem.
\newblock \emph{Doklady Mathematics}, 82(1), 540--542.
\newblock \doi{10.1134/S1064562410040101}.

\bibitem[{Leonov and Kuznetsov(2011)}]{LeonovK-2011-DAN}
Leonov, G. and Kuznetsov, N. (2011).
\newblock Algorithms for searching for hidden oscillations in the {A}izerman
  and {K}alman problems.
\newblock \emph{Doklady Mathematics}, 84(1), 475--481.
\newblock \doi{10.1134/S1064562411040120}.

\bibitem[{Leonov and Kuznetsov(2013)}]{LeonovK-2013-IJBC}
Leonov, G. and Kuznetsov, N. (2013).
\newblock Hidden attractors in dynamical systems. {F}rom hidden oscillations in
  {H}ilbert-{K}olmogorov, {A}izerman, and {K}alman problems to hidden chaotic
  attractors in {C}hua circuits.
\newblock \emph{International Journal of Bifurcation and Chaos}, 23(1).
\newblock \doi{10.1142/S0218127413300024}.
\newblock {a}rt. no. 1330002.

\bibitem[{Leonov and Kuznetsov(2014)}]{LeonovK-2014-book}
Leonov, G. and Kuznetsov, N. (2014).
\newblock \emph{Nonlinear Mathematical Models of Phase-Locked Loops. Stability
  and Oscillations}.
\newblock Cambridge Scientific Publisher.

\bibitem[{Leonov et~al.(2015{\natexlab{a}})Leonov, Kuznetsov, and
  Mokaev}]{LeonovKM-2015-EPJST}
Leonov, G., Kuznetsov, N., and Mokaev, T. (2015{\natexlab{a}}).
\newblock Homoclinic orbits, and self-excited and hidden attractors in a
  {L}orenz-like system describing convective fluid motion.
\newblock \emph{Eur. Phys. J. Special Topics}, 224(8), 1421--1458.
\newblock \doi{10.1140/epjst/e2015-02470-3}.

\bibitem[{Leonov et~al.(2012)Leonov, Kuznetsov, Yuldahsev, and
  Yuldashev}]{LeonovKYY-2012-TCASII}
Leonov, G., Kuznetsov, N., Yuldahsev, M., and Yuldashev, R. (2012).
\newblock Analytical method for computation of phase-detector characteristic.
\newblock \emph{IEEE Transactions on Circuits and Systems - II: Express
  Briefs}, 59(10), 633--647.
\newblock \doi{10.1109/TCSII.2012.2213362}.

\bibitem[{Leonov et~al.(2015{\natexlab{b}})Leonov, Kuznetsov, Yuldashev, and
  Yuldashev}]{LeonovKYY-2015-TCAS}
Leonov, G., Kuznetsov, N., Yuldashev, M., and Yuldashev, R.
  (2015{\natexlab{b}}).
\newblock Hold-in, pull-in, and lock-in ranges of {PLL} circuits: rigorous
  mathematical definitions and limitations of classical theory.
\newblock \emph{IEEE Transactions on Circuits and Systems--I: Regular Papers},
  62(10), 2454--2464.
\newblock \doi{10.1109/TCSI.2015.2476295}.

\bibitem[{Margaris(2004)}]{Margaris-2004}
Margaris, N. (2004).
\newblock \emph{Theory of the Non-Linear Analog Phase Locked Loop}.
\newblock Springer Verlag, New Jersey.

\bibitem[{Pliss(1958)}]{Pliss-1958}
Pliss, V.A. (1958).
\newblock \emph{Some Problems in the Theory of the Stability of Motion (in
  Russian)}.
\newblock Izd LGU, Leningrad.

\bibitem[{Shakhgil'dyan and Lyakhovkin(1966)}]{ShahgildyanL-1966}
Shakhgil'dyan, V. and Lyakhovkin, A. (1966).
\newblock \emph{Fazovaya avtopodstroika chastoty (in Russian)}.
\newblock Svyaz', Moscow.

\bibitem[{Shakhgil'dyan and Lyakhovkin(1972)}]{Shahgildyan-1972}
Shakhgil'dyan, V. and Lyakhovkin, A. (1972).
\newblock \emph{Sistemy fazovoi avtopodstroiki chastoty (in Russian)}.
\newblock Svyaz', Moscow.

\bibitem[{Shakhtarin(1969)}]{Shakhtarin-1969}
Shakhtarin, B. (1969).
\newblock Study of a piecewise-linear system of phase-locked frequency control.
\newblock \emph{Radiotechnica and electronika (in {R}ussian)}, (8), 1415--1424.

\bibitem[{Suarez et~al.(2012)Suarez, Fernandez, Ramirez, and
  Sancho}]{Suarez-2012}
Suarez, A., Fernandez, E., Ramirez, F., and Sancho, S. (2012).
\newblock Stability and bifurcation analysis of self-oscillating quasi-periodic
  regimes.
\newblock \emph{IEEE Transactions on microwave theory and techniques}, 60(3),
  528--541.

\bibitem[{Tricomi(1933)}]{Tricomi-1933}
Tricomi, F. (1933).
\newblock Integrazione di unequazione differenziale presentatasi in
  elettrotechnica.
\newblock \emph{Annali della R. Shcuola Normale Superiore di Pisa}, 2(2),
  1--20.

\bibitem[{Viterbi(1966)}]{Viterbi-1966}
Viterbi, A. (1966).
\newblock \emph{Principles of coherent communications}.
\newblock McGraw-Hill, New York.

\end{thebibliography}

\end{document}